\documentclass[12pt]{article}

\usepackage{amsmath,amssymb}
\usepackage{tikz}
\usepackage{color}
\usepackage[toc]{appendix}
\usepackage{graphicx}
\usepackage{fancyhdr}
\usepackage{enumitem}
\usepackage{bbm}

\usepackage{graphicx}
\usepackage{parskip}
\usepackage{float}
\usepackage{chngpage}
\usepackage{calc}
\usepackage{bigints}
\usepackage{tikz}
\usepackage{array}
\usepackage{bbm}
\usepackage{booktabs}
\usepackage{rotating}
\usepackage{multirow}
\usepackage{adjustbox}
\usepackage{tabularx}
\usepackage{enumitem}
\usepackage{verbatim}
\usepackage{mathtools}
\usepackage{ragged2e}
\usepackage[makeroom]{cancel}
\usepackage{enumitem}
\usepackage{caption}
\usepackage{mathtools}
\usepackage{bbm}
\usepackage{amsmath,amsfonts,amsthm,bm}
\usepackage{amsmath,amssymb}
\usepackage{appendix}
\usepackage{graphicx}
\usepackage{tikz}
\usepackage{verbatim}
\usepackage{amsthm}
\usepackage[english]{babel}
\usepackage{hyphenat}
\usepackage[makeindex]{imakeidx}
\usepackage{tikz}

\usetikzlibrary{datavisualization}
\usetikzlibrary{matrix}
\usetikzlibrary{datavisualization.formats.functions}

\newtheorem{theorem}{Theorem}[section]

\newtheorem{remark}[theorem]{Remark}
\newtheorem{assumption}[theorem]{Assumption}

\usepackage[colorlinks=true, citecolor=blue, urlcolor=blue]{hyperref}

\makeatletter
\def\section{\@startsection {section}{1}{\z@}{3.25ex plus 1ex minus
		.2ex}{1.5ex plus .2ex}{\large\bf}}
\def\subsection{\@startsection{subsection}{2}{\z@}{3.25ex plus 1ex minus
		.2ex}{1.5ex plus .2ex}{\normalsize\bf}}
\@addtoreset{equation}{section} 
\makeatother
\makeatletter
\newsavebox{\@brx}
\newcommand{\llangle}[1][]{\savebox{\@brx}{\(\m@th{#1\langle}\)}%
	\mathopen{\copy\@brx\kern-0.5\wd\@brx\usebox{\@brx}}}
\newcommand{\rrangle}[1][]{\savebox{\@brx}{\(\m@th{#1\rangle}\)}%
	\mathclose{\copy\@brx\kern-0.5\wd\@brx\usebox{\@brx}}}
\makeatother

\title{A small-time approximation of Girsanov's exponential and  short-time transition density estimations.}

\author{ Ramiro Scorolli\thanks{Dipartimento di Scienze Statistiche Paolo Fortunati, Università di Bologna, Bologna, Italy. \textbf{e-mail}: ramiro.scorolli2@unibo.it}}

\date{\today}

\begin{document}
	
	\maketitle
	
	\bigskip

	\begin{abstract}
	The main result of this article regards a small time approximation for the Girsanov's exponential. We prove that the latter is well described over
	short time intervals by the solution of a deterministic partial differential equation.The rate of convergence of the approximation is of order one in the length of the interval.
	As a possible application we show that our approximation can be used to obtain an estimation of the short-time transition density of a Langevin equation representing the dynamics of a Brownian particle moving under the influence of an external, non-linear force. Using this approach is equivalent to consider a random ordinary differential equation, where the dynamics of the particle is deterministic and given by the aforementioned force and the stochasticity enters through the initial condition.
	The result is in accordance with those obtained by discretization methods such as Euler-Maruyama, implicit Euler and  by approximating the associated Fokker-Planck-Smoluchowski equation.
	\end{abstract}
	
	
	
	\bigskip

	\section{Introduction and statement of the main results.}
	In this article we will derive a short-time approximation for the solution of the stochastic Cauchy problem:
	\begin{align}\label{SDE}
		\begin{cases}
			dM_t= f(\alpha+B_t)M_t dB_t,\;t\in [0, T ]\\
			M_0=1.
		\end{cases}
	\end{align}
	Here:
	\begin{itemize}
		\item $f:\mathbb R\to\mathbb R$ is a function satisfying assumption \ref{assumptions} (see below).
		\item $\{B_t\}_{t\in [0,T]}$ is a standard  $\mathbb R$-valued Brownian motion defined on some complete filtered probability space $(\Omega,\mathcal F,\{\mathcal F_t\},P)$ satisfying the usual conditions denotes an admissible filtration.
		\item $\alpha$ is a $\mathbb R$-valued, $\mathcal F_0$-measurable random variable with law given by $\mu$.
	\end{itemize}

	An straightforward application of It\^o formula \cite{oksendal2003stochastic} shows that the solution of (\ref{SDE}) is given by the following stochastic process:
	\small
	\begin{align}\label{Girsanov exponential}
		M_t:=\exp\bigg\{\int_0^t f(\alpha+B_s) dB_s-\frac{1}{2}\int_0^t|f(\alpha+B_s)|^2 ds\bigg\},\; t\in [0,T];
	\end{align}
	\normalsize
	the latter will be referred to as the \textit{Girsanov's exponential} .
	
	All throughout this article we will work under the following assumption:
	\begin{assumption}\label{assumptions}\quad
		\begin{itemize}
			\item $f\in C^2(\mathbb R)$;
			\item $f,f',f''$ are bounded;
			\item there exists a positive constant $\epsilon>0$ such that $f(x)>\epsilon, \forall x\in\mathbb R$.
		\end{itemize}
	\end{assumption}
	
	If the latter holds Girsanov's theorem (e.g. \cite{karatzas1998brownian}\cite{oksendal2003stochastic}) states that there exists a measure $Q$ on $(\Omega,\mathcal F)$ such that 
	\begin{align*}
		&W:[0,T]\times \Omega\to \mathbb R\\
		&(t,\omega)\mapsto B_t(\omega)-\int_0^t f(\alpha+B_s(\omega))ds 
	\end{align*}
	is a $Q$-Brownian motion and it holds that 
	\begin{align*}
		\frac{dQ}{dP}=M_T,
	\end{align*}
	i.e. the \textit{Girsanov's exponential} evaluated at the end-time $T$ is the Radon-Nikodym derivative of $Q$ with respect to $P$
	
	It is important to stress out the fact that assumption \ref{assumptions} will be needed in order to prove our main result and that Girsanov's theorem hold under much weaker conditions regarding $f$, in particular it's sufficient to assume that $f$ is bounded.
	
	Our approach shares has many similarities with those used in \cite{lanconelli2021small}  \cite{lanconelli2021wong} and \cite{hu1996wick} although the presence of a random diffusion coefficient force us to introduce certain modifications.
	
	The main novelty is that of connecting equation (\ref{SDE}) with a deterministic  partial differential equation that, just as the original equation, defines a Radon-Nikodym derivative. The latter will be done by exploiting the connection between the Wick product and It\^o integration and the fact that under certain conditions the Wick product behaves like a differential operator.
	
	Later on we will show that this approximation can be used in conjunction with Girsanov's theorem in order to describe the short-time transition density of an associated diffusion process that describes the behavior of a Brownian particle moving under the influence of an external force.
	
	We now present our main theorem;
	
	\begin{theorem}\label{main theorem}
		Let assumption \ref{assumptions} be in force, and for $ T >0$ let
		\begin{align*}
			[0, T ]\times \mathbb R\ni (t,x)\mapsto u(t,x)
		\end{align*}
		be the classical solution of the following deterministic Cauchy problem
		\begin{align}\label{pde}
			\begin{cases}
				\partial_t u=  -\partial_x\left[f\left(\alpha+x\right)u\right]+\frac{f\left(\alpha+x\right)u\,x}{ T },\;  t\in ]0, T ]\\
				u(0,x)=1.
			\end{cases}
		\end{align}
		
		Then there exists a measure $\mathcal Q$ on $(\Omega,\mathcal B)$ such that 
		\begin{align*}
			\mathcal M_{ T }(\omega):=u( T ,B_{ T }(\omega))=\frac{d\mathcal Q}{dP}
		\end{align*}
		i.e. $\mathcal M_{ T }$ it is the Radon-Nikodym derivative of $\mathcal Q$ with respect to $P$.
		
		Furthermore for any $p\geq 1$ it holds that 
		\begin{align}
			\big\|M_{ T }-\mathcal M_{ T }\big\|_{p}\leq CT .
		\end{align}
		where $\|\bullet\|_p$ denotes the $\mathbb L^p(\Omega)$-norm and $C$ is a suitable constant depending on $p$ and the other parameters of the model.
	\end{theorem}
	The proof and the motivation explaining why we consider a deterministic Cauchy problem will be postponed to the last section.

	\section{Estimations of the short-time transition densities.}
	
	Let's consider the following stochastic Cauchy problem
	\begin{align}\label{SDE2}
		\begin{cases}
			dX_t=f(X_t)dt+dB_t,\; t\in[0, T ]\\
			X_0=\alpha,
		\end{cases}
	\end{align}
	where again $f$ satisfies the assumption \ref{assumptions} and $\alpha$ is some $\mathcal F_0$ measurable random variable.
	This  could be seen as Langevin equation  representing  the dynamics of a Brownian particle influenced by an external force given by $f(\bullet)$ and where the initial state of the system is random. 
	
	If we let $f(x)=x$ and  $\alpha\equiv a\in \mathbb R$ in such a case the solution is given by the well-known Ornstein–Uhlenbeck process. The latter example is the only case, at least known to the author, in which the solution of an equation of the form of (\ref{SDE2}) can be obtained explicitly besides the trivial case where $f(\bullet)\equiv c\in\mathbb R$.
	
	\begin{remark}
		The fact that we are using a unit diffusion coefficient in (\ref{SDE2}) by no means entails a loss of generality in our discussion. 
		Consider a general diffusion given by: 
		\begin{align*}
			dY_t=b(Y_t)dt+\sigma(Y_t)dB_t,\; t\in [0,T].
		\end{align*}
		Under mild assumptions on $b$ and $\sigma$ we are able to transform the latter into a diffusion with the same structure as (\ref{SDE2}). The latter is sometimes referred to as  Lamperti’s Transformation ( e.g.  \cite{pavliotis2014stochastic}\cite{moller2010state}). Of course in our case, we would still have to check that the new drift coefficient satisfies the assumption \ref{assumptions}.
	\end{remark}
	
	In most cases we are more interested in the law of $X_t$ for any arbitrary $t\in [0,T]$ rather than the trajectory itself since the former gives a much more complete summary of the dynamics into play. In the context of Markovian stochastic processes (such as It\^o diffusions \cite{oksendal2003stochastic}) a complete description of the law of a given process can be obtained by means of the transition density and the Chapman-Kolmogorov equation \cite{schilling2014brownian} \cite{risken1996fokker}.
	
	For this reason in the following we will set ourselves to find approximations of the short-time transition density of the diffusion given by (\ref{SDE2}).

	In order to calculate the short-time transition density function $p(T,x|0,x')$ of an It\^o diffusion one in general can proceed in two different ways. 
	The first is to work on the stochastic differential equation, i.e. to employ approximation methods such like the Euler-Maruyama scheme and the Backward-Euler scheme (e.g. \cite{kloedenande1992platen}) on a single interval that would become arbitrarily small.
	In both cases one is naturally led to consider transformations of $B_T$, and from there it's easy to obtain a density function by inverting this transformation (see for instance \cite{horsthemke1975onsager}).
	
	On the other hand one could work instead with the associated Fokker-Planck-Smoluchowski equation. In \cite{haken1976generalized} \cite{risken1996fokker} the authors show how to obtain approximations of the transition probability by performing a discretization of the time interval. These estimations of the transition probabilities are then used to construct\textit{ path integrals solutions} of the Fokker-Planck equation.
	It's important to notice that in general the approximations of the transition probability obtained by the aforementioned methods are not actually densities, in particular they do not integrate to $1$.
	
	We propose a third alternative that makes use of Girsanov's theorem and the approximation obtained in the previous section.
	
	For the sake of completeness the estimations obtained using the aforementioned approaches will be presented:
	\subsection{Backward-Euler scheme}
	
	For a single interval the approximation is given by 
	\begin{align}\label{backward euler}
		X_{ T }=x'+f(X_{ T }) T +B_{ T }.
	\end{align}
	Notice that this describes a transformation $B_{ T }\mapsto X_{ T }$, whose inverse is given by
	
	\begin{align*}
		B_{ T }=X_{ T }-x'+f(X_{ T }) T.
	\end{align*}
	
	Then the density of the approximated solution is given by (e.g.\cite{casella2021statistical})
	\begin{align*}
		p( T ,x|0,x')=\frac{1}{\sqrt{2\pi T }}e^{-\frac{1}{2 T }[x-x'-f(x) T ]^2}\times|1-f'(x) T |.
	\end{align*}
	
	Using the fact that $|1-f'(x) T |\approx e^{-f'(x) T }$ we arrive to the following approximation of the transition density 
	
	\begin{align}\label{TD Backward Euler}
		p( T ,x|0,x')=\frac{1}{\sqrt{2\pi T }}e^{-\frac{1}{2 T }[x-x'-f(x) T ]^2-f'(x) T }.
	\end{align}

	\subsection{Euler-Maruyama scheme }
	Using the standard Euler-Maruyama method the approximation is given by
	
	\begin{align}\label{euler}
		X_{ T }=x'+f(x') T +B_{ T },
	\end{align}
	
	again this is a transformation $B_{ T }\mapsto X_{ T }$ and the inverse is given by 
	\begin{align*}
		X_{ T }-x'-f(x') T =B_{ T }.
	\end{align*}
	
	Then the approximated density is given by
	\begin{align}
		p( T ,x|0,x')=\frac{1}{\sqrt{2\pi T }}e^{-\frac{1}{2 T }[x-x'-f(x') T ]^2}.
	\end{align}

	\subsection{Using approximated Girsanov's exponential.}
	It's well known that Girsanov's theorem can be used to construct weak solutions of equations of the form of (\ref{SDE2}) (see for instance \cite{karatzas1998brownian} or \cite{oksendal2003stochastic}). The basic idea is that the law of the random variable $X_t(\omega)$ for some fixed time $t\in [0,T]$ equals the law of $B_t(\omega)$ under the measure  $dQ=M_TdP$, where again $M_T$ denotes the Girsanov's exponential evaluated at the end time.
	
	It's then reasonable to think that if we replace the  Girsanov's exponential $M_{ T }$ with our approximation $\mathcal M_T$ we could obtain an approximation (in law) $\mathcal X_T$ of $X_T$ as suggested by the following calculation
	\begin{align*}
		&|\mathbb E_P(g(X_T))-\mathbb E_P(g(\mathcal X_T))|\\
		&=|\mathbb E_P(g(\alpha+B_T)M_T)-\mathbb E_P(g(\alpha+B_T)\mathcal M_T)|\\
		&\leq \|g\|_{\infty}\|M_t-\mathcal M_T\|_1\\
		&\leq C T,
	\end{align*}
	for any   bounded Borel function $g:\mathbb R\to\mathbb R$.
	
	Again for an arbitrary bounded function $g:\mathbb R\to\mathbb R$ we have:
	\begin{align*}
		\mathbb E_P\left[g(\mathcal X_{ T })\right]&=\mathbb E_P\left[g(\alpha+B_{ T })\cdot\mathcal M_{ T }\right]\\
		&=\frac{1}{\sqrt{2\pi T }}\int_{\mathbb R}\int_{\mathbb R}g(a+b) \left(\frac{f\left(\Lambda^{-1}\left(\Lambda(a+b)- T \right)\right)}{f(a+b)}\right)\\
		&\quad\quad\dots\times e^{-\frac{1}{2 T }\left[\Lambda^{-1}\left(\Lambda(a+b)- T \right)-a\right]^2} db\; \mu(da),
	\end{align*}
	where we used the fact that $B_{ T }\sim N(0,T)$, $\mu$ denotes the law of the random variable $\alpha$, and we wrote $\mathcal M_T$ in explicit form (the reader is referred to  section  \ref{Proof} for the definition of $\Lambda$ and $\Lambda^{-1}$).
	
	This means that the approximated transition density function is given by
	\small
	\begin{align}\label{transition probability}
		\tilde{p}( T ,x|0,x')=\frac{1}{\sqrt{2\pi T }} \left(\frac{f\left(\Lambda^{-1}\left(\Lambda(x)- T \right)\right)}{f(x)}\right) e^{-\frac{1}{2 T }\left[\Lambda^{-1}\left(\Lambda(x)- T \right)-x'\right]^2}.
	\end{align}
	\normalsize
	By letting $y=\Lambda^{-1}\left(\Lambda(x)- T \right)$ (see equations (\ref{Lambda}) and (\ref{Lambda inverse})) it's straightforward to see that the latter is indeed a density with respect to the Lebesgue measure, i.e. is a non-negative function which integrates to $1$; this change of variables also entails the following.
	
	\begin{remark}
		The function $\tilde{p}$ is the probability density function of the random variable $\mathcal X_T:=\Lambda^{-1}(\Lambda(x'+B_T)+T)$. Furthermore notice that $\{\mathcal X_t\}_{t\in [0,T]}$ is the solution of the ordinary differential equation:
		\begin{align}\label{ODE}
			\begin{cases}
				\frac{d\mathcal X_t}{dt}=f(\mathcal X_t),\; t\in ]0,T]\\
				\mathcal X_0=x'+B_T,
			\end{cases}
		\end{align}
		which implies
		\begin{align}\label{our approach}
			\mathcal X_T=x'+B_T+\int_0^T f(\mathcal X_s)ds.
		\end{align}
		Then our approach is equivalent to solving an ODE (that describes the dynamics of the particle in a deterministic scenario) where the stochasticity enters through the initial condition (such equations are sometimes referred to as ``crypto-deterministic'' \cite{moyal1949stochastic}). It's important to stress out the fact that in certain cases \ref{ODE} can be solved under less restrictive conditions than assumption \ref{assumptions}, an straightforward example being $f(x)=x$.
		This approach lies somewhere in the middle between solving the actual SDE and using a discrete approximation since we only discretize the ``random component''( cf. (\ref{our approach}), (\ref{backward euler}) and (\ref{euler})).
		The latter can also be seen if one notices that (\ref{transition probability}) is the solution of 
		\begin{align}\label{liouville}
			\begin{cases}
				\partial_t \tilde{p}(t,x)=-\partial_x[f(x)\tilde{p}(t,x)],\; t\in [0,T]\\
				\tilde{p}(0,x)=\frac{1}{\sqrt{2\pi T}}e^{\frac{-(x-x')^2}{2T}},
			\end{cases}
		\end{align}
		evaluated at  $t=T$, i.e. a Liouville equation (e.g. \cite{soong1973random} \cite{banks2012uncertainty}) that describes the evolution of the probability density function of a process $\{X_t\}_{t\in [0,T]}$ that evolves in a deterministic setting (according to the force $f(\bullet)$) and  the initial condition is a Gaussian random variable with mean $x'$ and variance $T$.
		
		Equation (\ref{liouville}) can be compared to the Fokker-Planck-Smoluchowski equation associated with (\ref{SDE2}):
		\begin{align}\label{Fokker Planck}
			\begin{cases}
				\partial_t p(t,x)=-\partial_x[f(x)p(t,x)]+\frac{1}{2}\partial_{xx}^2 p(t,x),\; t\in[0,T]\\
				p(0,x)=\delta(x-x').
			\end{cases}
		\end{align}
		where $p(t,x)$ is a shorthand for $p(t,x|0,x')$.
		
		Notice that the first-order differential operator is the same in both equations, but the second-order operator associated with the heat-semigroup does not appear in (\ref{liouville}). Still the ``gaussianity'' enters in the equation through the initial condition.
		
	\end{remark}

	A solution for (\ref{Fokker Planck}) cannot be obtain in general but a class of equivalent approximations is available (e.g.\cite{haken1976generalized}). 
	Following \cite{risken1996fokker} one starts by considering a small time interval $[0,T]$ and hence: 
	\begin{align*}
		p(T,x|0,x')=[1+ T\mathcal L_{FP}+\mathcal O(T^2)]\delta(x-x')
	\end{align*}
	where $\mathcal L_{FP}\bullet=-\frac{\partial }{\partial x}f(x)\bullet+\frac{1}{2}\frac{\partial^2}{\partial x^2}\bullet$ .
	Using the integral representation of the delta function and integration by parts (see \cite{risken1996fokker}) one arrives to the following approximation 
	
	\begin{align}\label{Haken}
		p( T ,x|0,x')=\frac{1}{\sqrt{2\pi T }}e^{-\frac{1}{2 T }[x-\alpha-f(x) T ]^2-f'(x) T }.
	\end{align}
	It's worth to notice that the latter in general does not integrate to $1$, i.e. it's not an actual density.
	Furthermore notice that this approximation is non-unique as stressed out in \cite{risken1996fokker}.

	We can treat (\ref{liouville}) in an analogous way:
	\begin{align}\label{approx discrete}
		\tilde p(T,x|0,x')=[1- T\mathcal L+\mathcal O(T^2)]\tilde p(0,x),
	\end{align}
	where $\mathcal L \bullet=-\frac{\partial }{\partial x}f(x)\bullet$.
	
	The latter yields
	\begin{align*}
		\tilde p(T,x|0,x')=(2\pi T)^{-1/2}e^{-x^2/2T}[1-f'(x)T+f(x)x+\mathcal O(T^2)].
	\end{align*}
	
	Now notice that in \cite{haken1976generalized} the author argue that in order to obtain a ``correct'' solution one should have complete agreement up to and including terms of the order of $T$, same thing happens in \cite{wissel1979manifolds} where the author adds and subtract arbitrary terms of order $\mathcal O(T^2)$ to the analogous of (\ref{approx discrete}).
	
	Using the same logic we can add
	$f(x)^2\frac{T^2}{2}$ and neglect the $\mathcal O(T^2)$ inside the square brackets.
	Then, within the aforementioned level of accuracy, we have:
	\begin{align*}
		\tilde p(T,x|0,x')&=\frac{1}{\sqrt{2\pi T}}e^{-\frac{x^2}{2T}}\left[1-\left(f'(x)+\frac{f(x)x}{T}+\frac{f(x)^2T}{2}\right)T\right]\\
		&=\frac{1}{\sqrt{2\pi T}}e^{-\frac{1}{2T}[x-f(x)T]^2-f'(x)T},
	\end{align*}
	which coincides with (\ref{Haken}).
	These somewhat informal calculations allows us to see that the approximation (\ref{transition probability}) is part of the class of equivalent representations described by \cite{haken1976generalized}.
	
	\begin{remark}
		The fact that (\ref{Haken}) can be obtained by obtained by making further approximations to  (\ref{transition probability}) lead us to believe that ours is more precise, but since the exact solution isn't known we haven't been able to prove this conjecture rigorously. Further research will be done on this particular.
	\end{remark}
	
	\begin{remark}
		We can use our short-time transition density and Chapman-Kolmogorov equation in order to obtain path integral solutions (e.g.\cite{feynman2010quantum} \cite{haken1976generalized} \cite{wissel1979manifolds}) of (\ref{Fokker Planck}) for any arbitrary time $T>0$:
		\begin{align*}
			p(T,x|0,x')&=\lim_{\substack{N\to\infty\\N\tau=T}} (2\pi\tau)^{-N/2}\int_{\mathbb R^{N}}\Bigg(\left[\prod_{i=1}^N \left(\frac{f(\Lambda^{-1}(\Lambda(x_i)-\tau))}{f(x_i)}\right)dx_i\right]\\
			&\quad\quad\dots\times e^{-\frac{1}{2 \tau }\sum_{i=1}^{N}\left[\Lambda^{-1}\left(\Lambda(x_i)- \tau \right)-x_{i-1}\right]^2}\Bigg)
		\end{align*}
		where $x_0\equiv x'$.
	\end{remark}
	
	\begin{remark}\label{remark 2}
		If in (\ref{approx}) instead of replacing the Brownian motion appearing in the diffusion coefficient with $B_{ T }$ we replace it with $B_0=0$ we obtain a different partial differential equation, namely:
		\begin{align*}
			\begin{cases}
				\partial_t \tilde{u}=-f(\alpha)\partial_x\tilde{u}+\frac{f(\alpha)\tilde{u}x}{T}\\
				\tilde{u}(0,x)=1.
			\end{cases}
		\end{align*}
		
		Solving this equation and proceeding as explained in \ref{euristic} yields an alternative approximation of the Girsanov's exponential:
		\begin{align*}
			\tilde{\mathcal M}_{ T }(\omega)=e^{\frac{-1}{2 T }\left[B_{ T }(\omega)^2-(B_{ T }(\omega)-f(\alpha) T )^2\right]}.
		\end{align*}
		
		From there we have
		\begin{align*}
			\mathbb E_P[g(\tilde{\mathcal X}_{ T })]&=\frac{1}{\sqrt{2\pi T }}\int_{\mathbb R}\int_{\mathbb R} g(b+a) e^{\frac{-1}{2 T }[b-f(a) T ]^2}db\;\mu(da)\\
			&=\frac{1}{\sqrt{2\pi T }}\int_{\mathbb R}\int_{\mathbb R} g(b) e^{\frac{-1}{2 T }[b-a-f(a) T ]^2}db\;\mu(da)
		\end{align*}
		
		and hence the transition density is given by
		\begin{align*}
			p( T ,x|0,x')=\frac{1}{\sqrt{2\pi T }}e^{-\frac{1}{2 T }[x-x'-f(x') T ]^2},
		\end{align*}
		which is the density that  obtained by means of Euler-Mayurama method and coincides with equation $4.55$ of \cite{risken1996fokker}.
	\end{remark}
	
	\section{Proof of Theorem \ref{main theorem}}\label{Proof}
	
	We start this proof by explaining the motivations that led us to consider a deterministic Cauchy problem as an approximation for an SDE.
	
	\begin{remark}\label{euristic}
		Consider the stochastic differential equation (SDE) that defines the Girsanov's exponential, i.e.
		\begin{align*}
			\begin{cases}
				dM_t= f(\alpha+B_t)M_t dB_t,\;t\in [0,T]\\
				M_0=1.
			\end{cases}
		\end{align*}
		
		Using the connection the Wick product and the It\^o integration one could formally rewrite equation  (\ref{SDE}) as
		\begin{align}
			\begin{cases}
				\frac{dM_t}{dt}=f(B_t)M_t\diamond \frac{dB_t}{dt},\\
				X_0=1,
			\end{cases}
		\end{align}
		where $\diamond$ denotes the Wick product (e.g. \cite{simon2015p}) and $\frac{dB_t}{dt}$ denotes the (distributional) time-derivative of the Brownian motion (see for instance \cite{holden1996stochastic}\cite{kuo2018white}\cite{hida2013white} for a complete account on the connection between Wick product, white noise and It\^o integration).
		
		A natural approximation for  $\frac{dB_t}{dt}$ is given by  $\frac{B_{ T }}{ T }$ which amounts to consider the Wong-Zakai approximation \cite{wong1969riemann} with the coarser partition of the interval $[0, T ]$.
		Furthermore, following the same logic (in the spirit of the Backward-Euler scheme) we will replace the Brownian motion at time $t$ with the Brownian motion evaluated at the end of the interval, in remark \ref{remark 2} we will discuss what would happens if instead we replace $B_t$ with $B_0=0$ (in the spirit of the Euler-Mayurama method).

		With this in hand we can introduce the following ``approximating'' equation
		\begin{align}\label{approx}
			\begin{cases}
				\frac{d\mathcal {M}_t}{dt}= \left[f\left(\alpha+B_{ T }\right) \cdot \mathcal {M}_t\right]\diamond \left(\frac{B_{ T }}{ T }   \right), t\in ]0, T ]\\
				\mathcal M_0=1.
			\end{cases}
		\end{align}
		Under certain conditions the Wick product between a random variable $F$ and an element of the first Wiener chaos $I(g)$ is given by (e.g. \cite{holden1996stochastic},\cite{kuo2018white}) 
		\begin{align}
			F\diamond I(g)=F\cdot I(g)-D_g X
		\end{align}
		where $D_g$ stands for the directional Malliavin derivative in the direction of the deterministic function $g(\bullet)$ and $I(\bullet)$ denotes the Wiener-It\^o integral. Since $B_{ T }=\int_0^{ T }dB_s=I(1_{[0, T ]})$  we can formally write
		\small
		\begin{align*}
			\frac{d{\mathcal M_t}}{dt}= \frac{1}{ T }\bigg\{\left[f\left(\alpha+B_{ T }\right) \mathcal M_t\right]\cdot B_{ T }-\bigg\langle D\left[f\left(\alpha+B_{ T }\right) \mathcal M_t\right],1_{]0, T ]}\bigg\rangle \bigg\}.
		\end{align*}
		\normalsize
		If we now set ourselves to find a solution of the form:
		\begin{align*}
			\mathcal M(t,\omega)=u(t,B_{ T }(\omega)),\; t\in ]0, T ]
		\end{align*}
		for some continuous function $u:[0,T]\times \mathbb R\to\mathbb R$ to be determined, we can use the chain rule for Malliavin derivative to show that $u$ must solve
		\begin{align}\label{eq: pde deterministic}
			\begin{cases}
				\partial_t u=-\partial_x\left[f\left(\alpha+x\right)u\right]+ \frac{f\left(\alpha+x\right)u\,x}{ T },\; (t,x)\in ]0,T]\times\mathbb R\\
				u(0,x)=1.
			\end{cases}
		\end{align}
		It is important to stress the fact that this calculations are just formal and simply provide a candidate for the approximation; if we additionally assume that the random variable $F$ belongs to the Malliavin-Sobolev space $\mathbb D^{1,p},p>1$ (e.g  \cite{nualart2006malliavin}) then all the calculations are justified (see for instance \cite{holden1996stochastic} \cite{hu1996wick}).

	\end{remark}
	
	We are now ready to start the proof of the theorem;
	all throughout this proof we will denote with $ C$ a constant whose value may change at each appearance. 
	The starting point will be the following deterministic Cauchy problem (\ref{eq: pde deterministic})
	which is related with (\ref{SDE}) by the heuristic considerations discussed above.
	The transformation $v(t,x):=f(\alpha+x)u(t,x)e^{\frac{-x^2}{2 T }}$ and a simple application of the chain rule allows us to see that the following must hold
	\begin{align*}
		\begin{cases}
			\partial_t v= -f\left(\alpha+x\right)\partial_xv,\; t\in ]0, T ]\\
			v(0,x)=f(\alpha+x)e^{\frac{-x^2}{2 T }}=:\varphi(x).
		\end{cases}	
	\end{align*}
	
	This equation can be solved by a simple application of the method of characteristics which yields
	\begin{align*}
		v(t,x)=\varphi\left(\Lambda^{-1}\left(\Lambda(x)-t\right)\right),
	\end{align*}
	where $\Lambda:\mathbb{R}\to\mathbb R$ is given by
	\begin{align}\label{Lambda}
		\Lambda(x):=\int \frac{1}{f(\alpha+x)}dx
	\end{align}
	and $\Lambda^{-1}:\mathbb R\to\mathbb R$ satisfies
	\begin{align}\label{Lambda inverse}
		\Lambda^{-1}(\Lambda(x))=x;
	\end{align}
	notice that all this functions are well defined under assumption \ref{assumptions}.
	
	Inverting the transformation shows that the solution to (\ref{eq: pde deterministic}) is given by
	\begin{align}\label{RND}
		u(t,x)&=\left(\frac{f\left(\Lambda^{-1}\left(\Lambda(\alpha+x)-t\right)\right)}{f(\alpha+x)}\right)\times e^{-\frac{1}{2 T }\left[\left(\Lambda^{-1}\left(\Lambda(\alpha+x)-t\right)-\alpha\right)^2-x^2\right]}, \; t\in[0, T ],
	\end{align}
	as desired.
	
	We will now show that if we replace the spatial variable with the Brownian motion evaluated at time $T$ the resulting random variable $u( T ,B_{ T })(\omega)$ is a Radon-Nikodym derivative, i.e. it's a.s. positive and has mean $1$. 
	
	The positivity follows immediately from assumption \ref{assumptions} and the positivity of the exponential function.
	On the other hand using the fact that $B_{ T }\sim N(0, T )$ we have that
	
	\begin{align*}
		\mathbb E[u( T ,B_{ T })]&=\frac{1}{\sqrt{2\pi T }}\int_{\mathbb R}\int_{\mathbb R} \left(\frac{f\left(\Lambda^{-1}\left(\Lambda(a+b)-T\right)\right)}{f(a+b)}\right)\\
		&\quad\quad\dots\times e^{-\frac{1}{2 T }\left[\left(\Lambda^{-1}\left(\Lambda(a+b)-T\right)-a\right)^2-b^2\right]} e^{-\frac{b^2}{2 T }}db\;\mu(da),
	\end{align*} 
	
	changing variables under the prescription $y=\Lambda^{-1}\left(\Lambda(a+b)-T\right)$ yields
	
	\begin{align*}
		\mathbb E[u( T ,B_{ T })]=\frac{1}{\sqrt{2\pi T }}\int_{\mathbb R}\int_{\mathbb R}  e^{-\frac{(y-a)^2}{2 T }}dy\; \mu(da)=1.
	\end{align*}

	For the sake of simplicity and without loss of generality from now on we will assume that $\alpha\equiv 0$, since the general case can be obtained by straightforward modifications.
	
	We start by noticing that we can write
	\begin{align*}
		\mathcal M_{ T }&=	\exp\Bigg\{\left(\frac{\ln\left[f\left(\Lambda^{-1}\left(\Lambda(B_{ T })- T \right)\right)\right]-\ln\left[f(B_{ T })\right]}{ T }\right) T\\ &\quad\quad-\frac{1}{2 T }\left[\left(\Lambda^{-1}\left(\Lambda(B_{ T })- T \right)\right)^2-B_{ T }^2\right]\Bigg\},
	\end{align*}

	and then we and left to consider the following expression
	\begin{align*}
		\big\|\mathcal M_{ T }-M_{ T }\big\|_{p}&=\Bigg\|\exp\bigg\{\left(\frac{\ln\left[f\left(\Lambda^{-1}\left(\Lambda(B_{ T })- T \right)\right)\right]-\ln\left[f(B_{ T })\right]}{ T }\right) T\\ &\quad\quad-\frac{1}{2 T }\left[\left(\Lambda^{-1}\left(\Lambda(B_{ T })- T \right)\right)^2-B_{ T }^2\right]\bigg\}\\
		&\quad\quad-\exp\bigg\{\int_0^{ T } f(B_s)dB_s-\frac{1}{2}\int_0^{ T }|f(B_s)|^2 ds\bigg\}\Bigg\|_{p}.
	\end{align*}
	
	Using the elementary inequality\small $|e^X-e^Y|\leq (e^X+e^Y)\cdot|X-Y|$ \normalsize we have that
	\begin{align*}
		\|\mathcal M_{ T }-M_{ T }\|_{p}\leq \big\|(e^{\mathcal A}+e^{\mathcal B})\cdot |\mathcal A-\mathcal B|\big\|_{p}
	\end{align*}
	where 
	\begin{align*}
		\mathcal A&:=\left(\frac{\ln\left[f\left(\Lambda^{-1}\left(\Lambda(B_{ T })- T \right)\right)\right]-\ln\left[f(B_{ T })\right]}{ T }\right) T\\ 
		& \quad \quad-\frac{1}{2 T }\left[\left(\Lambda^{-1}\left(\Lambda(B_{ T })- T \right)\right)^2-B_{ T }^2\right],\\
		\mathcal B&:=\int_0^{ T } f(B_s)dB_s-\frac{1}{2}\int_0^{ T }|f(B_s)|^2 ds.
	\end{align*}
	
	An application of Hölder  and triangular inequality yields
	\begin{align}\label{inequality norm}
		\big\|\mathcal M_{ T }-M_{ T }\big\|_{p}&\leq \big\|\mathcal M_{ T }+M_{ T }\big\|_{p_1}\big\| \mathcal A-\mathcal B\big\|_{p_2}\nonumber\\
		&\leq \left(\big\|\mathcal M_{ T }\big\|_{p_1}+\big\| M_{ T }\big\|_{p_1}\right)\big\| \mathcal A-\mathcal B\big\|_{p_2},
	\end{align}
	where $\frac{1}{p_1}+\frac{1}{p_2}=\frac{1}{p}$.
	
	At this point we will need the following remark
	
	\begin{remark}\label{remark1}
		Under the assumption \ref{assumptions} it holds that for any $q\geq 1$  $\|\mathcal M_{ T }\|_{q}$ and $\|M_{ T }\|_{q}$ are bounded by a constant that does not change as $T\to 0$.
	\end{remark}
	
	\begin{proof}
		We will start by showing the latter for the Girsanov's exponential.
		
		Notice that
		\begin{align*}
			|M_{ T }|^q=\exp\bigg\{q\int_0^{ T } f(B_s)dB_s-\frac{q}{2}\int_0^{ T }|f(B_s)|^2 ds\bigg\},
		\end{align*}
		at this point we can define $g(\bullet):=qf(\bullet)$ and then the latter reads
		\begin{align*}
			|M_{ T }|^q=\exp\bigg\{\int_0^{ T } g(B_s)dB_s-\frac{1}{2q}\int_0^{ T }|g(B_s)|^2 ds\bigg\},
		\end{align*}
		or which is equivalent
		\begin{align*}
			|M_{ T }|^q&=\exp\bigg\{\int_0^{ T } g(B_s)dB_s-\frac{1}{2}\int_0^{ T }|g(B_s)|^2 ds\bigg\}\\	&\times\exp\bigg\{\frac{q-1}{2q}\int_0^{ T }|g(B_s)|^2 ds\bigg\}.
		\end{align*}
		At this point take the expectation of the expression above and notice that since $g$ is a bounded function, $\exp\big\{\frac{q-1}{2q}\int_0^{ T }|g(B_s)|^2 ds\big\}$ is bounded by some constant that remains unchanged as $T$ becomes arbitrarily small; the boundedness also implies that the Novikov's  condition \cite{karatzas1998brownian} holds and hence the first exponential has mean equal to $1$, this proves the result for $\|M_T\|_{q}$.
		
		We now show the boundedness of
		\small
		\begin{align*}
			\big\|\mathcal M_T\big\|_{q}&=	\Bigg\|\left(\frac{f\left(\Lambda^{-1}\left(\Lambda(B_{ T })- T \right)\right)}{f(B_{ T })}\right) e^{-\frac{1}{2 T }\left[\left(\Lambda^{-1}\left(\Lambda(B_{ T })- T \right)\right)^2-B_{ T }^2\right]}\Bigg\|_{q}.
		\end{align*}
		\normalsize
		
		Assumption \ref{assumptions} implies that the ratio is bounded by a constant that does not depend on $ T $, then 
		
		\begin{align*}
			\big\|\mathcal M_T\big\|_{q}&\leq C\bigg\| e^{-\frac{1}{2 T }\left[\left(\Lambda^{-1}\left(\Lambda(B_{ T })- T \right)\right)^2-B_{ T }^2\right]}\bigg\|_{q}\\
			&=C\left(\int_{\mathbb R} e^{-\frac{q}{2 T }\left[\left(\Lambda^{-1}\left(\Lambda(x)- T \right)\right)^2-x^2\right]} \frac{1}{\sqrt{2\pi T }}e^{-\frac{1}{2 T }x^2}dx\right)^{\frac 1 q}\\
			&=\left(\frac{C}{\sqrt{2\pi T }}\int_{\mathbb R} e^{-\frac{q}{2 T }\left[\left(\Lambda^{-1}\left(\Lambda(x)- T \right)\right)^2\right]}e^{\frac{q-1}{2 T }x^2}dx\right)^{\frac 1 q},
		\end{align*}
		where we used the fact that $B_{ T }\sim N(0, T )$.
		
		In order to evaluate the integral above we take a Taylor expansion of the function $\Lambda^{-1}(\bullet)$ around $x_0=\Lambda(x)$ which yields
		\begin{align}\label{taylor1}
			\Lambda^{-1}(x_0- T )&=x-f(x) T +\frac{1}{2}f(\xi)f'(\xi) T ^2,
		\end{align}
		for some $\xi\in ]x_0- T ,x_0[$, clearly such a $\xi$ depends on $x$ and to stress this dependence we write $\xi(x)$. Then plugging (\ref{taylor1}) into the last equality above yields 
		\begin{align}
			\big\|\mathcal M_T\big\|_{q}&\leq \left(\frac{C}{\sqrt{2\pi T }}\int_{\mathbb R} e^{-\frac{q}{2 T }\left[x-f(x) T +\frac{1}{2}f(\xi(x))f'(\xi(x)) T ^2\right]^2}e^{\frac{q-1}{2 T }x^2}dx\right)^{\frac{1}{q}}
		\end{align}
		straightforward calculations show that the right hand side above equals
		\begin{align*}
			\Bigg(\frac{C}{\sqrt{2\pi T }}&\int_{\mathbb R} \exp\bigg\{\frac{q}{2}\bigg[2xf(x)-xf(\xi(x))f'(\xi(x)) T -f(x)^2 T\\
			&+f(x)f(\xi(x))f'(\xi(x)) T ^2-\frac{1}{4}f(\xi(x))^2f'(\xi(x))^2 T ^3\bigg]\bigg\} e^{\frac{-x^2}{2 T }}dx\Bigg)^{\frac 1 q}.
		\end{align*}
		
		Now notice that the terms $-f(x)^2T -\frac{1}{4}f(\xi(x))^2f'(\xi(x))^2 T ^3$ are non-positive and can be upper-bounded by $0$ and that $f(x)f(\xi(x))f'(\xi(x)) T ^2$ is bounded and can be majored by some positive constant. We can further bound this expression by taking absolute values as follows
		\begin{align*}
			\big\|\mathcal M_T\big\|_{q}&\leq \left(\frac{C}{\sqrt{2\pi T }}\int_{\mathbb R} e^{\frac{q}{2}\left[2|x|f(x)+|x|f(\xi(x))|f'(\xi(x))| T \right]}e^{\frac{-x^2}{2 T }}dx\right)^{\frac 1 q}.
		\end{align*}
		
		Again, using the boundedness of $f$ and after reordering we obtain
		\begin{align*}
			\big\|\mathcal M_T\big\|_{q}&\leq \left(\frac{C}{\sqrt{2\pi T }}\int_{\mathbb R} e^{K|x|}e^{\frac{-x^2}{2 T }}dx\right)^{\frac{1}{q}}\\
			&=\left(C\left(\operatorname{erf}\left(\dfrac{K\sqrt{ T }}{\sqrt{2}}\right)+1\right)\mathrm{e}^\frac{K^2 T }{2}\right)^{\frac{1}{q}}\\
			&\leq \left(2C e^{K^2/2}\right)^{\frac{1}{q}},
		\end{align*}
		where we used the boundedness of the error function notice that the constants $C$ and $K$ do not change as  $T$ becomes arbitrarily small.
		
	\end{proof}
	
	We can now come back to the main proof, using remark \ref{remark1} it would suffices to prove the result for
	\begin{align}\label{A-B}
		\|\mathcal A-\mathcal B\|_{p_2}=&\Bigg\|\left(\frac{\ln\left[f\left(\Lambda^{-1}\left(\Lambda(B_{ T })- T \right)\right)\right]-\ln\left[f(B_{ T })\right]}{ T }\right) T\nonumber\\ &-\frac{1}{2 T }\left[\left(\Lambda^{-1}\left(\Lambda(B_{ T })- T \right)\right)^2-B_{ T }^2\right]\nonumber\\
		&-\int_0^{ T } f(B_s)dB_s+\frac{1}{2}\int_0^{ T }|f(B_s)|^2 ds\Bigg\|_{p_2}.
	\end{align}

	Notice that by the Fundamental theorem of calculus  and assumption \ref{assumptions} one has:
	
	\begin{align}\label{first term}
		&\left(\frac{\ln\left[f\left(\Lambda^{-1}\left(\Lambda(B_{ T })- T \right)\right)\right]-\ln\left[f(B_{ T })\right]}{ T }\right)\nonumber\\
		&=-\frac{d}{dx}\ln\left[f(\Lambda^{-1}(x))\right]\big|_{x=\Lambda(B_{ T })}+\mathcal O( T ^2)\nonumber\\
		&=-f'(B_{ T })+\mathcal O( T ^2).
	\end{align}
	
	Now using (\ref{taylor1}) and letting $x=B_{ T }(\omega)$ yields
	\small
	\begin{align*}
		\left(\Lambda^{-1}\left(\Lambda(B_{ T })- T \right)\right)^2=&\bigg[ B_{ T }^2 -2B_{ T }f(B_{ T }) T +B_{ T }\mathcal O(T^2) +f(B_T)T^2+\mathcal O(T^3)\bigg].
	\end{align*}
	
	Then the following equality holds
	\begin{align}\label{second term}
		&-\frac{1}{2 T }\left[\left(\Lambda^{-1}\left(\Lambda(B_{ T })- T \right)\right)^2-B_{ T }^2\right]\nonumber\\
		&=f(B_{ T })B_{ T }-B_{ T }\mathcal O( T )-\frac{1}{2}f(B_{ T })^2 T -\mathcal O( T ^2)\nonumber\\
		&=f(B_{ T })\diamond B_{ T }+f'(B_{ T }) T -B_{ T }\mathcal O( T )-\frac{1}{2}f(B_{ T })^2 T +\mathcal O(T^2),
	\end{align}
	
	where we  used the fact that $f(B_{ T })\cdot B_{ T }=f(B_{ T })\diamond B_{ T }+\big\langle D f(B_{ T }),1_{[0, T ]}\big\rangle$.
	
	Plugging (\ref{first term}) and (\ref{second term}) into (\ref{A-B}) we obtain 
	\begin{align*}
		\|\mathcal A-\mathcal B\|_{p_2}
		&\leq \bigg\|f(B_{ T })\diamond B_{ T }-\int_0^{ T } f(B_s)dB_s\\
		&\quad \quad+\frac{1}{2}f(B_{ T })^2  T -\int_0^{ T }|f(B_s)|^2ds -B_T\mathcal O(T)+\mathcal O(T^2)\bigg\|_{p_2}.
	\end{align*}
	and using the triangular inequality yields
	\begin{align*}
		\|\mathcal A-\mathcal B\|_{p_2}
		&\leq \bigg\|f(B_{ T })\diamond B_{ T }-\int_0^{ T } f(B_s)dB_s\bigg\|_{p_2}\\
		&+\frac{1}{2}\bigg\|f(B_{ T })^2  T -\int_0^{ T }|f(B_s)|^2ds \bigg\|_{p_2}+\mathcal O( T ^{3/2}).
	\end{align*}
	
	Now notice that 
	\begin{align*}
		\bigg|f(B_{ T })^2 T -\int_0^{ T }|f(B_s)|^2ds \bigg|&=\bigg|\int_0^{ T }f(B_{ T })^2-f(B_s)^2ds \bigg|\\
		&\leq \int_0^{ T }|f(B_{ T })^2-f(B_s)^2|ds\\
		&\leq \int_0^{ T }C\cdot|B_{ T }-B_s|ds
	\end{align*}
	
	where we used the fact that $f(\bullet)^2$ is a Lipschitz continuous function.
	By taking the $\mathbb L^{p_2}$-norm and applying Minkowski inequality for integrals we obtain
	\begin{align}\label{lebesgue integral}
		\frac{1}{2}\bigg\|f(B_{ T })^2 T -\int_0^{ T }|f(B_s)|^2ds \bigg\|_{p_2}&\leq C\int_0^{ T }\|B_{ T }-B_s\|_{p_2}ds\nonumber\\
		&=C \int_0^{ T } ( T -t)^{1/2}dt\nonumber\\
		&\leq C T ^{3/2}.
	\end{align}
	
	On the other hand we have that 
	\begin{align*}
		\bigg|f(B_{ T })\diamond B_{ T }-\int_0^{ T } f(B_s)dB_s\bigg|&=\bigg|f(B_{ T })\diamond \int_0^{ T }dB_s-\int_0^{ T } f(B_s)dB_s\bigg|\\
		&=\bigg|\int_0^{ T }f(B_{ T })\delta B_s-\int_0^{ T } f(B_s)dB_s\bigg|\\
		&=\bigg|\int_0^{ T }f(B_{ T })-f(B_s)\delta B_s\bigg|\\
		&=|\delta(u)|,
	\end{align*}
	where  $u:=((t,\omega)\ni [0,T]\times \Omega \mapsto f(B_{ T }(\omega))-f(B_t(\omega)))$, and $\delta$ denotes the Skorohod integral \cite{nualart2006malliavin}\cite{janson1997gaussian}.
	
	Let  $H:=\mathbb L^2([0,\mathcal T])$, for some fixed $\mathcal T\geq T$, then proposition $1.5.8$ of \cite{nualart2006malliavin} states that for $u\in\mathbb D^{1,p}(H), p>1$ it holds that:
	\begin{align*}
		\|\delta(u)\|_{p}\leq c_p\left(\|\mathbb E(u)\|_{H}+\|Du\|_{\mathbb L^p(\Omega:H\otimes H)}\right).
	\end{align*}
	In our case this reads:
	\begin{align*}
		&\bigg\|f(B_{ T })\diamond B_{ T }-\int_0^{ T } f(B_s)dB_s\bigg\|_{p_2}\\
		&\leq C \Bigg(\|\mathbb E(u{1}_{[0,T]})\|_H\\
		&+\mathbb{E}\left[\left(\int_{ [0,\mathcal T]^2 }\big|D_s\left(u_t{1}_{[0,T]}(t)\right)\big|^2dsdt\right)^{\frac{p_2}{2}}\right]^{\frac{1}{p_2}}\Bigg),
	\end{align*}
	
	where the Malliavin derivative is given by
	\begin{align*}
		D_s\left(u_t{1}_{[0,T]}(t)\right)&=1_{[0,T]}(t)\left[f'(B_{ T })1_{[0,T]}(s)-f'(B_t)1_{[0,t]}(s)\right],
	\end{align*}
	for $(t,s)\in [0,\mathcal T]^2$, from where the membership of $u$ to $\mathbb D^{1,p}(H)$ is clear by the boundedness of $u$ and its Malliavin derivative.
	
	It follows that 
	\begin{align}\label{Ito integral}
		&\bigg\|f(B_{ T })\diamond B_{ T }-\int_0^{ T } f(B_s)dB_s\bigg\|_{p}\nonumber\\
		&\leq C\left(\left(\int_0^T \mathbb E\left[|f(B_T)-f(B_s)|^2\right]ds\right)^{\frac{1}{2}}+\mathcal O(T)\right)\nonumber\\
		&\leq CT,
	\end{align}
	where used the Lipschitz continuity of $f$ and the boundedness of the Malliavin derivative.
	
	The following toy-example seems to suggest that this bound is sharp, even though the function under consideration does not satisfy the assumption \ref{assumptions}
	
	\begin{remark}
		Letting $f(x)=x$ and $p=2$ we have that 
		\begin{align*}
			\bigg\|B_{ T }\diamond B_{ T }-\int_0^{ T } B_sdB_s\bigg\|_{\mathbb L^2(\Omega)}&=\bigg\|B_{ T }^{\diamond 2}-\frac{1}{2}B_{ T }^{\diamond 2}\bigg\|_{\mathbb L^2(\Omega)}\\
			&=\frac{1}{2}\|B_{ T }^{\diamond 2}\|_{\mathbb L^2(\Omega)}\\
			&=\frac{1}{2}\bigg\|\int_0^{ T }B_sdB_s\bigg\|_{\mathbb L^2}\\
			&=\frac{1}{2}\left[\int_0^{ T }sds\right]^{1/2}\\
			&=C T .
		\end{align*}
		suggesting that our bound is as sharp as it can possibly be.
	\end{remark}
	
	Finally from (\ref{inequality norm}), using (\ref{lebesgue integral}), (\ref{Ito integral}) and remark \ref{remark1} we conclude that
	\begin{align}
		\|\mathcal M_{ T }-	M_{ T }\|_{p}\leq  C T ,
	\end{align}
	as desired.

	\section{Discussion}
	The Girsanov's theorem is a powerful and quite useful tool in stochastic analysis. In this article we showed that the Girsanov's exponential can be approximated over short time intervals by the solution of a deterministic partial differential equation. This approximation can be used in order to estimate the transition density of the solution of a Langevin equation for small time intervals.
	
	Even though our estimations involve inverting arbitrary functions which in general doesn't lead to closed forms, our approach could be used in conjunction with interpolation methods in order to obtain numerical simulations. Further research will be done in this direction.
	Finally, the extension of this result to the multi-dimensional case is a work in progress.
	
	\section*{Acknowledgement}
	I would like to thank my supervisor Prof. Alberto Lanconelli for proposing me to work on this argument and for carefully going through all the proofs.
	I would also thank Prof. Alberto Scorolli, Prof. Fernando Tohmé and Simona Pace for the helpful discussions and for carefully reading the manuscript.
	
	\bibliographystyle{IEEEtran}
	\bibliography{ref}

\end{document}